\documentclass[10pt]{article}
\usepackage[english]{babel}
\usepackage{amsmath,amsthm}
\usepackage{amsfonts}
\usepackage[none]{hyphenat}
\usepackage{xcolor}
\usepackage{graphicx}
\usepackage[font={small}]{caption}
\usepackage[explicit]{titlesec}
\usepackage{anyfontsize}
\usepackage{colortbl}
\titleformat{\section}{\Large}{\textbf{\thesection .}}{1em}{\textbf{#1}}

\addto{\captionsenglish}{}
\newtheorem{thm}{Theorem}[section]

\newtheorem{con}[thm]{Conjecture}

\theoremstyle{definition}

\theoremstyle{remark}

\numberwithin{equation}{section}

\begin{document}
\title{Finding Prime Numbers\linebreak as Fixed Points of Sequences}%
\author{\vspace{-3,5\baselineskip}Enrique Navarrete$^{1}$\\
Daniel Orellana$^{2}$}\footnotetext[1]{Universidad Pontificia Bolivariana; enrique.navarrete@upb.edu.co}\footnotetext[2]{Interijento Consulting, daniel.orellana@interijento.net}%
\renewcommand*{\thefootnote}{\arabic{footnote}}
\makeatletter
\def\maketitle{%
\bgroup
\par\centering{\LARGE\@title}\\[3em]%
\ \par
{\@author}\\[4em]%
\egroup
}

\maketitle
\begin{center}
{\vspace{-2em} \scriptsize July 20, 2019 \vspace{-2em}}
\end{center}

\begin{abstract}
  In this note we describe a method for finding prime numbers as fixed points of particularly simple sequences. Some basic calculations show that success rates for identifying primes this way are over 99.9\%. In particular, it seems that the set of odd primes can be obtained as fixed points of the sequence which we call $A(1)$, the sequence of smallest divisors of triangular numbers, where the divisors are positive numbers that have not yet appeared in the sequence.\\[1em]
\end{abstract}
\section{Introduction}\label{sec1}

Throughout history, there have been many attempts to find either formulas or polynomials that produce prime numbers.  The most notable recent results of this type are due to Matiyasevich (1971) and Jones, Sato, Wada and Wiens (1976).  The first result, Matiyasevich's Theorem, is an existence result which states that there exists a finite set of Diophantine equations from which the set of prime numbers can be obtained.  The second result constructs an explicit polynomial of degree 25 in 26 variables whose positive values is the set of prime numbers.  Since there is plenty of bibliography on these topics, the interested reader is encouraged to look at the references.

A different recent approach to produce primes has been developed by Rowland (2008), who found a recurrence that produces primes and the value 1 (both with repetitions).  In this note we also detour from closed-form prime-generating approaches (\emph{ie}. formulas, polynomials), to find a class of sequences where odd primes are produced as fixed points of these sequences.  We describe the method below and some computational results in the next sections.

\section{Key sequences and Fixed-Point Algorithm}\label{sec2}

The method we use to generate odd primes in the correct positions in the sequence of positive numbers relies on the construction of a class of sequences which we will denote by $\mathcal{A}(p)$, where sequences $A(p)$ in the class are themselves indexed by the odd primes $p$. A sequence $A(p)$ can be described as the sequence of smallest divisors of $p$-multiples of triangular numbers, where the divisors are positive numbers that have not yet appeared in the sequence.

 We describe the method as follows:
\begin{enumerate}
  \item Choose an odd prime $p$, produce its multiples $kp$ and take the partial sums
  \[q(n)=\sum^{n-1}_{k=0}(kp)\text{, }n \geq 1; \text{ these are the $p$-multiples of triangular numbers,}\]
  OEIS A000217.
  \item Define a sequence with offset 1 as follows: let $a(1) = 1$ and let $a(n)$ be the smallest positive number not yet in the sequence such that $a(n)$ divides $q(n)$.   Denote the sequence constructed this way by $A(p)$.  Note that this sequence is well defined, that $a(1) = 1$ divides $q(1) = 0$, and that $a(2) = q(2) = p$.
  \item Under some qualifications mentioned below, we get odd prime numbers as fixed points of the sequence $A(p)$, \emph{ie}. as the terms $a(n)$ such that  $a(n) = n$.
\end{enumerate}
To illustrate, Table \ref{table:tab1} shows the first few terms of A(7).

\begin{table}[h!]
  \centering
  \small\begin{tabular}{ccccc}
    Fixed Points \ \ \ \ \ \ & \hspace{-0.2cm}$n$ & Mult of 7 & Sum $q(n)$ & $a(n)$ \\
  \end{tabular}
  \begin{tabular}{p{2cm}p{0.8cm}p{1cm}p{1.2cm}|c|}
    \cline{5-5}
    fixed point & 1 & 0 & 0 &  \ 1 \ \\
     & 2 & 7 & 7 & 7 \\
    fixed point & 3 & 14 & 21 & 3 \\
     & 4 & 21 & 42 & 2 \\
    fixed point & 5 & 28 & 70 & 5 \\
     & 6 & 35 & 105 & 15 \\
     & 7 & 42 & 147 & 21 \\
     & 8 & 49 & 196 & 4 \\
     & 9 & 56 & 252 & 6 \\
     & 10 & 63 & 315 & 9 \\
    fixed point & 11 & 70 & 385 & 11 \\
     & 12 & 77 & 462 & 14 \\
     fixed point & 13 & 84 & 546 & 13 \\
     & 14 & 91 & 637 & 49 \\
     & 15 & 98 & 735 & 35 \\
     & 16 & 105 & 840 & 8 \\
    fixed point & 17 & 112 & 952 & 17 \\
     & 18 & 119 & 1071 & 51 \\
    fixed point & 19 & 126 & 1197 & 19 \\
     & 20 & 133 & 1330 & 10 \\
     & 21 & 140 & 1470 & 30 \\
     & 22 & 147 & 1617 & 33 \\
    fixed point & 23 & 154 & 1771 & 23 \\
     & 24 & 161 & 1932 & 12 \\
     & 25 & 168 & 2100 & 20 \\
    \cline{5-5}
  \end{tabular}
  \caption{The first terms of $A(7)$.}\label{table:tab1}
\end{table}
As we can see from the table, the fixed points of the sequence are the first values of the set of odd primes (with the exception of 7, which does not appear as a fixed point since $a(2) = p$ in all sequences $A(p)$,  as mentioned above).\footnote[3]{Note that the value 1 also appears as a fixed point but is not a prime; however, it will not be counted as a \textquotedblleft false negative\textquotedblright\  or Type II error in Section \ref{sec4} since, by definition, $a(1) = 1$ in sequences $A(p)$.}

Even though this method is indeed very precise for finding primes, as well as highly efficient, it does have a minimal amount of error, \emph{ie}. the fixed points of the sequences $A(p)$ miss some primes.  Fortunately, as will be discussed below, the primes missed by $A(p)$ can be detected by some other sequence $A(p')$,  $p' \neq p$.  Besides this minimal amount of error, the only pitfall encountered in the method is that it also produces nonprimes as fixed points (which we call \textquotedblleft false negatives\textquotedblright\ from now on, using statistical terminology). To separate both kind of errors (miss primes $vs$. detect nonprimes), in Section \ref{sec3} we describe results leaving out false negatives, which we will return to in Section \ref{sec4}. Figure \ref{fig:Fig1} summarizes the two types of mistakes using terms from statistics or binary classification in general.\footnote[4]{Assuming the null hypothesis \textquotedblleft H$_0$: the number $n$ is prime\textquotedblright, a Type I error occurs when H$_0$ is rejected when it is true, or in our method, when a prime number $n$ is not detected as a fixed point of the sequence $A(p)$. A Type II error occurs when H$_0$ is false and is not rejected; that is, when a number is nonprime and is detected as a fixed point.  In hypothesis testing, a positive result corresponds to rejecting the null hypothesis, while a negative result corresponds to failing to reject it; the term \textquotedblleft false\textquotedblright\ means the conclusion inferred is incorrect. Hence a Type I error is also lnown as a \textquotedblleft false positive\textquotedblright\ and a Type II error as a \textquotedblleft false negative\textquotedblright.}

\begin{figure}[!h]
  \centering
\footnotesize\begin{tabular}{m{2.5cm}c}
   & \ \ \ \ PRIME  \ \ \ \ \ \ \ \ \ \ NONPRIME \\
  \ \ DETECT \linebreak\mbox{\ \ \ \ DON'T DETECT}&\normalsize\begin{tabular}{|c|c|}

                               \hline
                               Correct&False negative\\
                               \hline
                               False positive&Correct\\
                               \hline
                             \end{tabular}

\end{tabular}
  \caption{Classification Matrix.}\label{fig:Fig1}
\end{figure}


\section{Additional Observations and Results for Correctly Identified Primes}\label{sec3}
To check success rates of sequences $A(p$) to detect primes, the first 10,000 terms of such sequences were computed for $p = $ 3, 5, 7, 11, 13, 17, 19, 23, 41, 97 and 199. Table \ref{table:tab2} below shows that for $p$ as low as $p = 97$,  the rate of success for detecting primes is as large as 99.92\%.  Before describing the results, we highlight an additional observation that appeared in all cases computed, which we state as a conjecture.

\begin{con}
The fixed points of sequences $A(p)$, $p$ an odd prime, are odd.\footnote[5]{Actually, this also seems to be the case for $A(1)$, the sequence of smallest divisors of triangular numbers not yet in the sequence, as will be discussed in Section \ref{sec5}.}
\end{con}
Another observation stated as a conjecture is the following:

\begin{con}
A prime $n$ in the sequence $A(p)$ appears either as the fixed point    $n = a(n)$ or as $n = a(n+1)$.
\end{con}

For example, the earliest failure in all cases computed occurs in $A(3)$ for the prime 17, which appears at $a(18)$.

Note that if Conjecture 3.2 is true, it implies that if we allow a prime $n$ in the sequence $A(p)$ to appear either as the term $a(n)$ or as the term $a(n+1)$ (not too restrictive as compared to other methods), then the success rate for detecting primes is 100\%.\footnote[6]{Keeping again in mind the issue of false negatives, which will be discussed in Section \ref{sec4}.}
	
\begin{figure}[h!]
  \centering
  \includegraphics[width=12cm]{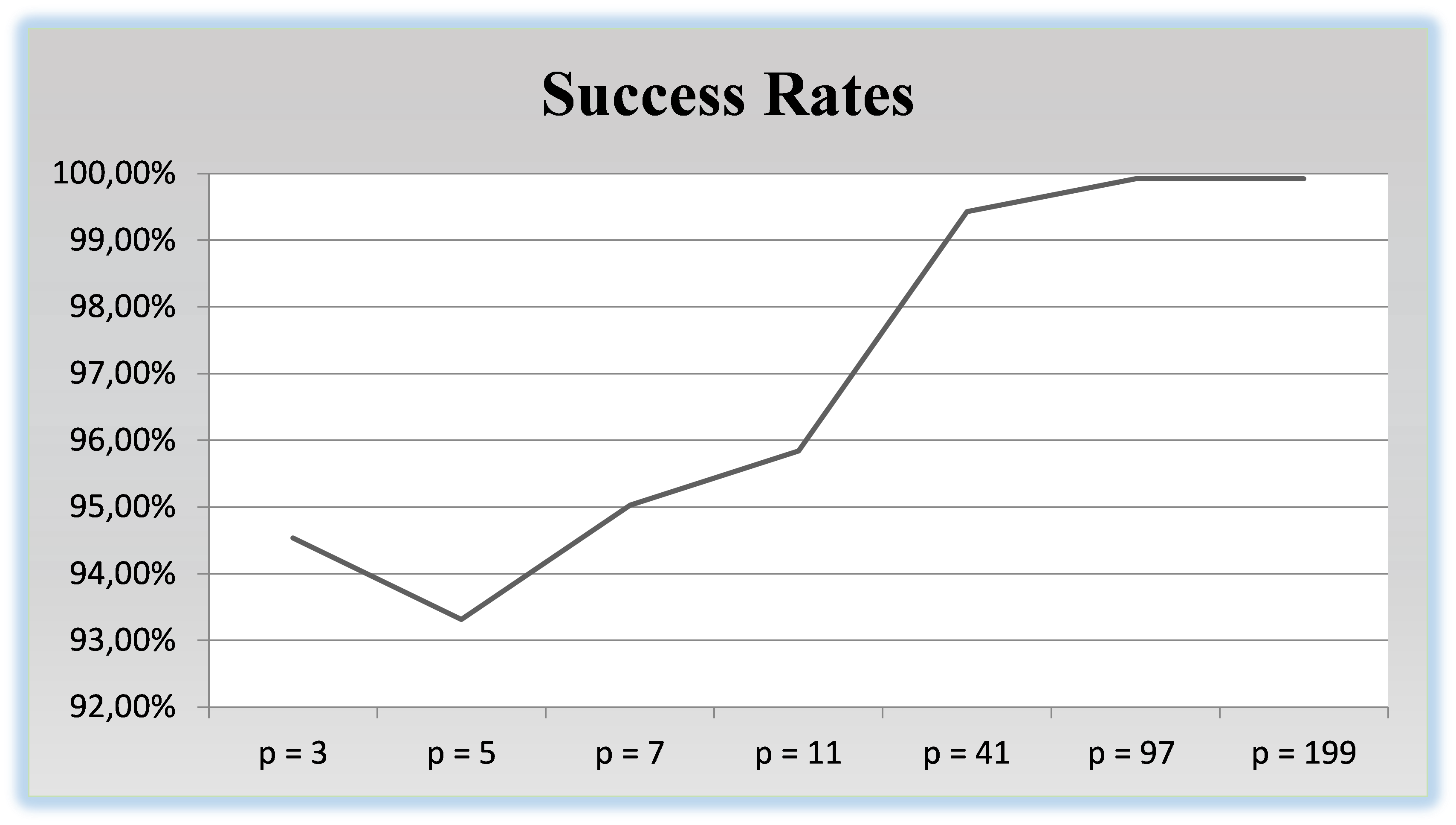}\\
  \caption{Success rates in sequences $A(p)$.}\label{fig:fig2}
\end{figure}

If we stick to the original requirement that prime $n$ must appear in the term $a(n)$ of the sequence $A(p)$,\footnote[7]{That is, the prime of value $n$ in its correct place in the sequence of positive numbers (not the $n^{\text{th}}$ prime).} then the natural question to ask is how well the sequences $A(p)$ perform in detecting prime numbers depending upon the selection of $p$.  As Table \ref{table:tab2} and Figure \ref{fig:fig2} show, the success rates are very favorable and even though they are initially high (94.5\% in $A(3)$), they seem to keep increasing with the selection of $p$ in $A(p)$.

\begin{table}[h!]
\hspace{-2cm}\small\begin{tabular}{|l|c|c|c|c|c|c|c|}
  \hline
  Success Rates of Sequences $A(p)$ & $p=3$ & $p=5$ & $p=7$ & $p=11$ & $p=41$ & $p=97$ & $p=199$ \\
  \hline
     A) Number of Matches $n=a(n)$ & 1160 & 1145 & 1166 & 1176 & 1220 & 1226 & 1226 \\
  \hline
     B) Number of Matches $n=a(n+1)$ & 67 & 82 & 61 & 51 & 7 & 1 & 1\\
  \hline
    C) Total Primes in 10,000 terms$^{*}$ & 1227 & 1227 & 1227 & 1227 & 1227 & 1227 & 1227 \\
    \hline
    Success Rates (A$/$C) & 94.54\% & 93.32\% & 95.03\% & 95.84\% & 99.43\% & 99.92\% & 99.92\%\\
  \hline
  \multicolumn{8}{l}{\scriptsize$ ^{\ast}$ Excluding 2 and $p$}
\end{tabular}
  \caption{Succes rates of sequences $A(p)$.}\label{table:tab2}
\end{table}

Note that there are 1229 primes in the first 10,000 terms of the positive numbers; however, in computing success rates for detecting odd primes in the sequences $A(p)$ we only consider 1227 primes since we exclude 2 as it is even and we exclude $p$ itself since it appears as the second term of every sequence $A(p)$, as discussed in Section \ref{sec1} (hence it is never a fixed point).

As Table \ref{table:tab2} shows, success rates for identifying odd primes by the first 10,000 terms of the sequences $A(p)$ start at 95.5\% for $p = 3$ and reach up to 99.9\% for $p = 199$.  The table distinguishes between primes identified by the correct term of the sequence, \emph{ie}. \textquotedblleft matches $n = a(n)$\textquotedblright\  and primes detected by the following term, $a(n+1)$.  Naturally, success rates in the table are computed using \textquotedblleft true\textquotedblright matches $n = a(n)$, \emph{ie}. the fixed points, since as mentioned above, using both types of matches would give a success rate of 100\% for all sequences $A(p)$, as can be checked from the table.

Note from Table \ref{table:tab2} that, aside from the issue of false negatives, the sequences A(97) and A(199) make only one mistake in identifying correctly the first 1227 primes (excluding 2 and $p$).  It turns out that the prime they both miss is the 844$^{\text{th}}$ prime, which is 6529.  The natural question that arises here is whether the same primes are missed by the sequences $A(p)$, as well as what kind of primes are missed. \emph{ie}.  primes of the form $4k+1$, $4k-1$, etc.  These and other questions will be pursued further on.


\section{Results taking into account False Negatives}\label{sec4}

We obtain the results below when false negatives are considered, \emph{ie}. nonprimes detected as fixed points of sequences $A(p)$.

\begin{table}[h!]
\hspace{-2.2cm}\small\begin{tabular}{|l|c|c|c|c|c|c|c|}
  \hline
  False Negatives in Sequences $A(p)$ & $p=3$ & $p=5$ & $p=7$ & $p=11$ & $p=41$ & $p=97$ & $p=199$ \\
  \hline
     A) Number of false negatives in 10,000 terms & 1179 & 1233 & 1248 & 1415 & 1478 & 1518 & 1526 \\
  \hline
    B) Number of nonprimes in first 10,000 terms & 8771 & 8771 & 8771 & 8771 & 8771 & 8771 & 8771 \\
    \hline
    \% (A$/$B) & 13.44\% & 14.06\% & 14.23\% & 16.13\% & 16.85\% & 17.31\% & 17.40\%  \\
  \hline
\end{tabular}

  \caption{Number of False Negatives in sequences $A(p)$.}\label{table:tab3}
\end{table}

As we see from Table \ref{table:tab3}, the percentage of false negatives in the first 10,000 terms of $A(p)$ increases with $p$, with percentages ranging from 13.4\% to 17.4\%.  These percentages are the ratio of the number of nonprimes detected as fixed points in the first 10,000 terms of $A(p)$ over the total number of nonprimes in the first 10,000 positive numbers.\footnote[8]{Recall that the value 1 appears as a fixed point in all sequences $A(p)$ but is not counted as a false negative since it is by definition that $a(1) = 1$ in $A(p)$ (see footnote 3).}

Another way of summarizing the performance of the sequences $A(p)$ in detecting primes is by means of the classification matrix, an example of which is shown below for the first 10,000 terms of $A(3)$.  Naturally, we would like the errors (off-diagonals) to be as small as possible.

\begin{figure}[!h]
  \centering
\footnotesize\begin{tabular}{m{3cm}c}
   & \ \ \ \ PRIME  \ \ \ \ \ \ \ \ \ \ NONPRIME \\
    \ \ DETECT\ \ \linebreak\mbox{\ \ \ \ DON'T DETECT}&\footnotesize
    \begin{tabular}{|c|c|}
                               \hline
                               \ \ \ \ \ 94.54\%\ \ \ \ \ &\ \ \ \ \ 13.44\%\ \ \ \ \ \\
                               \hline
                               \ \ \ \ \ 5.46\%\ \ \ \ \ &\ \ \ \ \ 86.56\%\ \ \ \ \ \\
                               \hline
                             \end{tabular}

\end{tabular}
  \caption{Classification Matrix for $A(3)$.}\label{fig:Fig3}
\end{figure}

From Tables \ref{table:tab2} and \ref{table:tab3} we see there is a tradeoff between success rates and percentage of false negatives since both increase with $p$.

Even though the percentages of false negatives seem to be high, we see from the actual sequences computed that many of these wrongly identified nonprimes are multiples of 3 or 5, so they are easily detected. It would be interesting to compute the percentage of false negatives after removing these nonprimes; it is likely that these percentages will decrease considerably.


\section{Results for nonprime values of of $p$ in $A(p)$}\label{sec5}

What kind of sequences $A(p)$ do we get when we select $p$ nonprime? Table \ref{table:tab4} below provides a brief summary.

\begin{table}[h!]
  \centering
    \small
  \begin{tabular}{|c|c|c|c|}

    \hline
    $p$ & $q(n)$ & OEIS & $A(p)$ \\
    \hline
    1 & Triangular Numbers & A000217 & A111273$^{\ast}$ \\
    \hline
    2 & Oblong numbers & A002378 & Positive Integers \\
    \hline
    4 & 4 times the triangular numbers & A046092 & Positive Integers \\
    \hline
    6 & 6 times the triangular numbers & A028896 & Positive Integers \\
    \hline
    9 & 9 times the triangular numbers & A027468 & Same as for $p=3$ \\
    \hline
    \multicolumn{4}{l}{\scriptsize$ ^{\ast} A(1)$ is A111273 if 0, the first triangular number, is omitted.}
  \end{tabular}
  \caption{Some sequences $A(p)$ for $p$ nonprime.}\label{table:tab4}
\end{table}

Table \ref{table:tab4} shows sequences $q(n)$, their classification number in OEIS, and the corresponding sequences $A(p)$ for some values of $p$ nonprime. We see that we don't get any new sequences $A(p)$ except for $p = 1$.\footnote[9]{Actually, for some values of  $p$ which are product of primes we get sequences $A(p)$ for some $p$ prime but with the first terms permuted.} In this case the sequence $q(n)$ is given by the triangular numbers  $n(n+1)/2 = \{0, 1, 3, 6, 10, 15, 21, \ldots\}$ (A000217 in OEIS). To get the sequence $A(1)$, if we drop 0, the first triangular number, then $A(1)$ is A111273 in OEIS, which starts as $\{1, 3, 2, 5, 15, 7, 4, 6, 9, 11, $ $22, 13, \ldots\}$. The fixed points of this sequence start with $\{1, 9, 25, 49, 57, 65, \ldots\}$. This sequence of fixed points, which is A113659 in OEIS, seems to produce some odd square numbers. However, if we compute $A(1)$ including the first triangular number 0, then by selecting 1 as the divisor of 0 (or $a(1) = a(2) = 1$), we get $A(1) = \{1, 1, 3, $ $2, 5, 15, 7, \ldots\}$, which is A111273 with an initial value of 1 appended. This shifts A111273 to the right by 1 and now the fixed points of the sequence are the odd prime numbers as in our previous results. Interestingly, the success rate of the shifted sequence up to the first 10,000 terms is 100\% since in this case all 1228 odd primes in the first 10,000 positive numbers are detected as fixed points. On the other hand, the number of false negatives is 1532, for a rate of 17.47\%, the largest one in our results.\footnote[10]{Note that considering the first triangular number 0 when computing $A(1)$ is consistent with our method described in Section \ref{sec2}, which computes the partial sum $q(n)$ of multiples $kp$ starting from $k = 0$.}

After checking the first 10,000 terms of the sequence $A(1)$, we state the following conjecture:

\begin{con}
Let $A(1)$ be the sequence of smallest divisors of triangular numbers not yet in the sequence with $a(1) = a(2) = 1$. Then the set of odd primes can be obtained as fixed points of this sequence.\footnote[11]{We don't use the term \textquotedblleft the fixed points\textquotedblright\  due to the false negatives.}
\end{con}

\section{Conclusion and Further Lines of Research}\label{sec6}

The method described for finding odd prime numbers is simple and efficient since, unlike other very intricate algorithms or formulas, it only requires sums of multiples of primes, divisions to obtain the terms in the sequences $A(p)$, and an algorithm to check that the smallest divisor $a(n)$ of $q(n)$ is a positive number not yet in the sequence.

With respect to precision, Section \ref{sec3} shows that success rates for finding primes are very high for the first 10,000 terms of the sequences $A(p)$. Furthermore, we have seen that these rates improve as the selection of $p$ in $A(p)$ increases, to the extent of missing just one prime in the first 10,000 terms of the sequences A(97) and A(199). Hence, we would expect that no primes will be missed after a sufficiently high selection of $p$. However, until further results are obtained, we state these expectations as the following conjecture:

\begin{con}
Let $p$ be a fixed odd prime and $a(n)\in A(p)$ the smallest positive number not yet in the sequence such that $a(n)$ divides $\displaystyle{q(n)=\sum^{n-1}_{k=0}(kp)}$,
$a(1) = 1$.  Then for large $p$, the set of odd primes $($excluding $p)$ can be obtained\\[1em]
as fixed points of the sequence $A(p)$.\footnote[12]{We exclude the value $p$ since recall that  $p = a(2)$ in all sequences $A(p)$.  Also, we don't use the term \textquotedblleft the fixed points\textquotedblright\  due to the false negatives.}$^{,}$\footnote[13]{Actually, at the time of submission, all 1227 primes (excluding 2 and $p$) were correctly identified by the first 10,000 terms of $A(p)$ for $p = 541$, the 100$^{th}$ prime, for a success rate of 100\%.}
\end{con}

Regarding the errors present in the method, the \textquotedblleft true errors\textquotedblright, or primes not detected, seem to be manageable, since the computations show that if some prime is not picked up by a particular sequence $A(p)$, it is detected by another (and there are plenty of sequences $A(p)$ to work with). More important seems to be the \textquotedblleft noise\textquotedblright\  introduced by false negatives, or nonprimes detected as fixed points, and even though they may be a drawback, many of them are multiples of 3 or 5, as mentioned previously, so they can be easily detected and discarded.\footnote[14]{Note that if Conjecture 3.1 is true, false negatives are odd.}

It seems that further work should concentrate on determining what kind of primes are missed by the sequences $A(p)$, as well as what are the \textquotedblleft best\textquotedblright\  individual sequences or combinations of sequences to minimize the percentage of false negatives.

\newpage
\section*{References}

{\small

[1] R. L. Goodstein and C. P. Wormell, Formulae for primes, \emph{The Mathematical} \hspace*{0.5cm}\emph{Gazette} \textbf{51} (1967) 35 - 38.

[2] James Jones, Formula for the $n$th prime number, \emph{Canadian Mathematical Bulletin} \hspace*{0.5cm}\textbf{18} (1975) 433 - 434.

[3] James Jones, Daihachiro Sato, Hideo Wada, and Douglas Wiens, Diophantine \hspace*{0.5cm}representation of the set of prime numbers, \emph{The American Mathematical Monthly} \hspace*{0.5cm}\textbf{83} (1976) 449 - 464.

[4] Yuri Matiyasevich, Diophantine representation of the set of prime numbers (in \hspace*{0.5cm}Russian), \emph{Doklady Akademii Nauk SSSR} \textbf{196} (1971).

[5] Rowland, Eric S. (2008), \textquotedblleft A Natural Prime-Generating Recurrence\textquotedblright, \emph{Journal of} \hspace*{0.5cm}\emph{Integer Sequences}, \textbf{11}: 08.2.8, arXiv:0710.3217.

[6] N. J. A. Sloane, The On-Line Encyclopedia of Integer Sequences, published \linebreak \hspace*{0.5cm}electronically at http://oeis.org.

[7] C. P. Willans, On formulae for the $n$th prime number, \emph{The Mathematical Gazette} \hspace*{0.5cm}\textbf{48} (1964) 413  - 415.

}



\end{document}